\input amstex
\magnification=\magstep1
\input amssym.def
\input amssym
\baselineskip 14 pt
\def \pop#1{\vskip#1 \baselineskip}

\font\gr=cmbx12

\def \et{\operatorname {et}}

\def \Im{\operatorname {Im}}
\def \Pic{\operatorname {Pic}}
\def \det{\operatorname {det}}
\def \ker{\operatorname {ker}}
\def \coker{\operatorname {coker}}
\def \Spec{\operatorname {Spec}}

\def \Sp{\operatorname {Sp}}
\def \ker{\operatorname {ker}}
\def \new{\operatorname {new}}
\def \geom{\operatorname {geom}}
\def \ProfGrp{\operatorname {ProfGrp}}\pop {1}
\def \ss{\operatorname {ss}}
\def \n{\operatorname {n}}
\def \Sets{\operatorname {Sets}}
\def \min{\operatorname {min}}
\def \deg{\operatorname {deg}}
\def \char{\operatorname {char}}
\def \Card{\operatorname {Card}}

\par
\centerline
{\bf \gr On Complete families of curves with a given fundamental}
\centerline
{\bf \gr group in positive characteristic}

{\pop 2}
\noindent
\centerline {\bf \gr Mohamed Sa\"\i di}

{\pop 2}
\noindent

\centerline {\bf \gr Abstratct}
{\pop 2}
\par
We prove in this paper, that complete families of smooth, and projective 
curves, of genus $g\ge 2$, in characteristic $p>0$, with a constant geometric 
fundamental group, are isotrivial.
\pop {2}
\par
\noindent
{\bf \gr 0. Introduction.}\rm\ Let $k$ be an algebraically closed field, and
let $X$ be a complete, irreducible, and smooth curve over $k$, of genus $g$.
The structure of the \'etale fundamental group $\pi_1(X)$ of $X$ 
is well understood, if $\char (k)=0$, thanks to the 
Riemann existence Theorem. Namely, it is isomorphic to the profinite 
completion $\Gamma_g$, of the topological fundamental group of a compact, 
orientable, topological surface, of genus $g$. In particular, the structure
of $\pi_1(X)$ depends only on $g$, in this case. 
In the case, where $\char (k)=p>0$, 
the structure of the full $\pi_1(X)$ is far from being understood. However,
we understand, in this case, the structure of some quotients of $\pi_1(X)$. 
Assume $\char (k)=p>0$. Let $\pi _1^p(X)$ (resp. 
$\pi _1(X)^{p'}$) be the maximal pro-p-quotient of $\pi_1(X)$ 
(resp. its maximal prime-to-p quotient). The following results are well known:

\par (1)\ The fundamental group $\pi_1(X)$, in characteristic $p>0$, 
is a quotient of the group $\Gamma_g$. In particular, $\pi_1(X)$ is 
topologically finitely generated.

\par (2)\ The structure of $\pi _1(X)^{p'}$ is well known, by Grothendieck's 
specialization theory for fundamental groups (cf. [SGA-1], X). Namely, 
it is isomorphic to the maximal prime-to-$p$ quotient of $\Gamma_g$.

\par (3)\ The structure of $\pi _1^p(X)$ is well known, by Shafarevich 
theorem (cf. [Sh]). Namely, it is a free pro-p-group on $r:=r_X$ generators, 
where $r_X$ is the $p$-rank of the curve $X$.

\par
Apart from these results, very little is known about the structure 
of the (geometric) fundamental group of curves, in positive characteristic.
\par
The anabelian geometry (or philosophy), as initiated by Grothendieck 
(cf. [Gr]), predicted that the structure of the arithmetic fundamental group, 
of hyperbolic curves over number fields, should depend on the isomorphy 
type of the curve in discussion. It came as a surprise when Tamagawa proved, 
in [Ta-3], such an anabelian statement, for hyperbolic affine curves, 
defined over a finite field of characteristic $p>0$. 
Tamagawa's result suggests, that anabelian phenomena may even hold, 
for the geometric fundamental group of complete curves, 
over arbitrary algebraically closed fields of characteristic $p>0$, 
which would explain to some extend, the complexity of
$\pi_1$ in positive characteristic. 
In this paper we give a new evidence to these expectations.

\par
In order to get an idea about the complexity of $\pi _1$ of proper curves,
in positive characteristic, we introduce the notion of fundamental group,
for points in the moduli space of curves. Let $\Cal M_g\to \Bbb F_p$ be 
the coarse moduli space of proper and smooth curves, of genus $g\ge 2$, 
in charactersitic $p>0$. 
It is well known that $\Cal M_g$ is a quasi-projective,
and geometrically irreducible variety. Let $L$ be an algebraically closed
field, of characteristic $p$. Then $\Cal M_g(L)$ is 
the set of isomorphism classes
of irreducible, proper, and smooth curves of genus $g$, over $L$. 
For a point $\overline x\in \Cal M_g(L)$, let 
$C_{\overline x}\to \Spec L$ be a curve classified by $\overline x$, and let
$x\in \Cal M_g$ be a point, such that $\overline x:\Spec L\to \Cal M_g$, 
factors through $x$. We
define the geometric fundamental group $\pi _1(x):=\pi _1(C_{\overline x})$, 
of the point $x$, as the fundamental group of the curve $C_{\overline x}$.
We remark that the structure of $\pi _1(x)$, as a profinite group, depends 
only on $x$, and not on the concrete geometric point 
$\overline x\in \Cal M_g(L)$ used to define it. 
Indeed, if $\overline {k(x)}$ is 
an algebraic closure of the residue field $k(x)$ at $x$, 
and $C_x$ is a curve classified by
$\Spec \overline {k(x)}\to \Cal M_g$, then 
$C_{\overline x}\simeq C_x\times _{\overline {k(x)}} L$ is the base change of 
$C_x$ to $L$. Hence, $\pi _1(C_{\overline x})\simeq \pi _1(C_x)$
by the geometric invariance of the fundamental group for proper varieties.

\pop {.5}
\par

A key tool in the study of $\pi_1$, is Grothendieck's specialization 
theory for fundamental groups (cf. [SGA], X). Let $y\in \Cal M_g$ be a point, 
which specializes in $x\in \Cal M_g$. Then, by Grothendieck's specialization 
theorem, there exists a surjective, continuous, homomorphism 
$\Sp_{y,x}:\pi _1(y)\to \pi_1(x)$. Concerning this specialization homomorphism, we have the following result:

\pop {.5}
\par
\noindent
{\bf \gr Theorem (Sa\"\i di, Pop, Raynaud, Tamagawa):}\rm\ {\sl 
Let $\overline {\Bbb F}_p$ be an algebraic closure of the finite field 
$\Bbb F_p$. Let $x\in \Cal M_g\times _{\Bbb F_p}\overline {\Bbb F}_p$ 
be a {\bf closed} point, and let $y\in \Cal M_g\times _{\Bbb F_p}
\overline {\Bbb F}_p$ be a point which specializes in $x$. Then, 
the specialization homomorphism $\Sp_{y,x}:\pi _1(y)\to \pi_1(x)$ 
is not an isomorphism.}

\pop {.5}
\par 
The above theorem was proven by Pop and Sa\"\i di, in [Po-Sa], in the
special case where the point $x$ corresponds to a curve, 
having an absolutely simple jacobian, with $p$-rank equal to $g$ or $g-1$,
and by Raynaud, in [Ra-1], in the case $g=2$, and the case of supersingular 
curves of arbitrary genus $g>2$ (i.e. curves whose jacobian is isogenous 
to a product of supersingular elliptic curves), 
and finally by Tamagawa, in [Ta-1], in the general case.

\pop {.5}
\par 
This result suggests that the structure, of the geometric fundamental 
group $\pi_1$, is far from being constant on the moduli space 
$\Cal M_g$, in characteristic $p>0$, 
much contrary to the characteritic $0$ case.
Let $k$ be an algebraically closed field, of characteristic $p>0$. 
Let $\Cal S\subset \Cal M_g\times _{\Bbb F_p}k$ be a $k$-subvariety. We say that
the geometric fundamental group $\pi_1$ is
{\it constant} on $\Cal S$, if for any two points $x$ and $y$ of $\Cal S$, such
that $y$ specializes in $x$, the corresponding specialization 
homomorphism $\Sp_{y,x}:\pi _1(y)\to \pi_1(x)$ is an isomorphism. This, 
in particular, would imply that all points of $\Cal S$ have isomorphic 
geometric fundamental groups. We say that $\pi_1$ is {\it not constant} 
on $\Cal S$, if the contrary holds, namely: there exists two points $x$ and 
$y$ of $\Cal S$, such that $y$ specializes in $x$, and such that the 
corresponding specialization homomorphism 
$\Sp_{y,x}:\pi _1(y)\to \pi_1(x)$, is not an isomorphism. 
The above Theorem implies, in particular, that in the case 
$k=\overline {\Bbb F}_p$, the moduli space 
$\Cal M_g\times _{\Bbb F_p}\overline {\Bbb F}_p$ does not contain, positive
dimensional $\overline {\Bbb F}_p$-subvarieties, 
on which $\pi_1$ is constant. It is thus natural to ask the following question:

\pop {.5}
\par
\noindent
{\bf \gr Question 6.3.}\rm\ Let $k$ be an algebraically closed field of 
characteristic $p>0$. Does $\Cal M_g\times _{\Bbb F_p}k$ contain 
$k$-subvarieties, of positive dimension $>0$, on which $\pi_1$ is constant?

\pop {.5}
\par

Our main result answering the above question is the following:

\pop {.5}
\par
\noindent
{\bf \gr Theorem (6.4).}\rm\ {\sl Let  $k$ 
be an algebraically closed field, of characteristic $p$. 
Let $\Cal S\subset \Cal M_g\times _{\Bbb F_p}k$ 
be a {\bf complete} $k$-subvariety of $\Cal M_g\times _{\Bbb F_p}k$. Then, the 
fundamental group $\pi_1$ is {\it not constant} on $\Cal S$.}

\pop {.5}
\par
Note, that it is well known that $\Cal M_g\times _{\Bbb F_p}k$ contains 
complete subvarieties (cf. [Oo-1], and [Fa-Lo], for example). 
In the case where, the generic points of the subvariety $\Cal S$, are 
contained in the locus of ordinary curves, i.e. curves having maximal 
$p$-rank equal to $g$, this result is well known, and due to Szpiro, 
Raynaud, and Moret-Bailly (cf. [Sz], and [Mo]). It can also be deduced
from Oort's result, on complete families of abelian varieties of dimension $g$,
with constant $p$-rank equal to $g-1$ (cf. [Oo], 6.2), 
if the generic points of the subvariety $\Cal S$, correspond 
to curves with $p$-rank equal to $g-1$.
One may 
wonder, whether there exists non-isotrivial complete smooth families of curves, of genus $g$, 
with constant $p$-rank. For otherwise, this would directly imply the above 
theorem, since the $p$-rank is encoded in the isomorphy type of the 
fundamental group. It turns out that such families exist. In [Oo-1], Oort 
constructed an example of a non-isotrivial complete smooth family of curves, 
of genus $3$, having constant $p$-rank equal to $0$. In section 5, we extend 
Oort's argument, in order to construct such examples, for any genus $g\ge 2$.

\pop {.5}
\par
For the proof of Theorem 6.4, it is easy to reduce to the case where
$\Cal S$ is a complete curve. In this case we prove, the following more 
precise result:

\pop {.5}
\par
\noindent
{\bf \gr Theorem 6.6.}\rm\ {\sl Let $k$ be an algebraically closed field, of 
characteristic $p>0$. Let $S$ be a smooth, complete, and irreducible $k$-curve.
Let $f:X\to S$ be a {\bf non-isotrivial}, proper and smooth family of curves, of 
genus $g\ge 2$. Then, there exists a finite \'etale cover 
$S'\to S$ of $S$, an \'etale cover $Y'\to X':=X\times _SS'$, of 
degree prime to $p$, and
a closed point $s_0\in S'$, such that
the $p$-rank of the geometric fibre $Y'_{k(\bar s_0)}\to k(\bar s_0)$, of $Y'$
above the point $s_0$, is strictly smaller, that the $p$-rank
of the generic geometric fibre  $Y'_{k(\bar \eta)}\to k(\bar \eta)$ of $Y'$
above the generic point $\eta$ of $S'$.}

\pop {.5}
\par
The main ingredients we use, in order to prove Theorem 6.6, are: first, 
Raynaud's theory of theta divisors in positive characteristic. 
Secondly, the Theorem of Szpiro, Raynaud, and Moret-Bailly, 
on the isotriviality of complete families of ordinary abelian varieties. And finally, 
a recent result of Tamagawa, on the equi-characteristic deformation 
of generalized Prym varieties. Finally, the statement of Theorem 6.4
can be easily generalized to the case where we consider the (geometric) tame 
fundamental group, see Theorem 6.10.
\pop {.5}
\par
This paper is organized as follows. In sections 1, and 2, 
we review Raynaud's theory of theta divisors in characteristic $p>0$, 
and its application to the study of the $p$-rank of cyclic, 
of order prime-to-p, \'etale covers of curves. 
In section 3, we explain Tamagawa's result
on  the equi-characteristic deformation of generalized Prym varieties. 
In section 4, we recall the results of Szpiro, Raynaud, Moret-Bailly, 
on complete families of abelian varieties with 
constant maximal $p$-rank. In section 5, we extend an argument of Oort,
in order to construct a complete family of smooth curves, for every genus 
$g\ge 2$, which has constant $p$-rank equal to $0$. In section 6, we prove our 
main result, using the results exposed in sections 1, 2, 3, and 4.
\pop {.5}
\par
First, I would like, very much, to thank F. Oort for answering my questions
concerning the content of  5. I would like also to thank M. Raynaud and 
A. Tamagawa for their comments on an earlier version of this paper. Finally, 
I would like to thank F. Pop for many interesting discussions we had on the 
subject, and for his interest in this work.

\par

\pop {.5}
\par
\noindent
{\bf \gr 1. The sheaf of locally exact differentials in characteristic 
$\bold {p>0}$ and its theta divisor.}

\pop {.5}
\par
In this section we review, mainly following Raynaud, the definition of the 
sheaf of locally exact differentials associated to a smooth algebraic curve in 
positive characteristic and its theta divisor
(cf. [Ra], 4, and [Ta], 1, for futher generalisations). 
Let $X$ be a proper smooth and connected algebraic curve 
of genus $g_X:=g\ge 2$, over an algebraically closed field $k$ of 
characteristic $p>0$. Consider the 
following cartesian diagram:

$$
\CD
X^1 @ >>> X \\
@VVV      @VVV  \\
\Spec k @>F>> \Spec k \\
\endCD
$$

\pop {.5}
\par
\noindent
where $F$ denotes the absolute Frobenius
morphism. 
The projection $X^1\to X$ is a 
scheme isomorphism. In particular, $X^1$ is a 
smooth and proper curve of genus $g$. The absolute Frobenius morphism
$F:X\to X$ induces in a canonical way a morphism  
$\pi:X\to X^1$ called the {\it relative Frobenius}, which is a radicial 
morphism of degree $p$. The canonical differential 
$\pi _*d:\pi _*\Cal O_X\to 
\pi _*\Cal \Omega ^1_X$ is a morphism of $\Cal O_{X^1}$-modules, its image
$B_X:=B:=\Im (\pi _*d)$ is the {\it sheaf of locally exact differentials}. 
One has the following exact sequence:
$$0\to \Cal O_{X^1}\to \pi _*\Cal O_X\to B\to 0$$

\noindent
and $B$ is a vector bundle on $X^1$ of rank $p-1$. 

\par
Consider the {\it Cartier operation} $c : \pi _*(\Cal 
\Omega ^1_X)\to \Cal \Omega ^1_{X^1}$,
 which is 
a morphism of $\Cal O_{X^1}$-modules. The kernel $\ker (c)$ of $c$ is equal 
to $B$, and the following sequence of $\Cal O_{X^1}$-modules is exact
(cf. [Se], 10):
$$0\to B\to \pi _{*}(\Cal \Omega^1_{X})\to \Cal \Omega ^1_{X^1}\to 0$$ 

\par
Let $\Cal L$ be a {\it universal Poincar\'e bundle} on $X^1\times _k J^1$, 
where $J^1:=\Pic ^0(X^1)$ is the Jacobian variety of $X^1$. The restriction of 
$\Cal L$ to $X^1\times \{a\}$, for any $a\in J^1(k)$, is
isomorphic to the degree zero line bundle $\Cal L_a$ which is the image of $a$
under the natural isomorphism $J^1(k)\simeq \Pic^0(X^1)$ ($\Cal L$ is 
normalized in such a way that $\Cal L_{0_{J^1}}\simeq \Cal O_{X^1}$). 
Let $h:X^1\times _k
J^1\to X^1$, and $f:X^1\times _kJ^1\to J^1$, be the canonical projections. As
$R^if_*(h^*B\otimes L)=0$ for $i\ge 2$, the total direct image $Rf_*(h^*B
\otimes L)$, of $(h^*B\otimes L)$ by $f$, can be realized by 
a complex $u:\Cal M^0\to \Cal M^1$ of length
$1$, where $\Cal M^0$ and $\Cal M^1$ are vector bundles on $J^1$, $\ker u=
R^0f_*(h^*B\otimes L)$, and $\coker u=R^1f_*(h^*B\otimes L)$. Moreover,
as the Euler-Poincar\'e characteritic $\chi (h^*B\otimes L)=0$, the vector 
bundles $\Cal M^0$ and $\Cal M^1$ have the same rank. In [Ra], Th\'eor\`eme
4.1.1, Raynaud proved the following theorem:

\pop {.5}
\par
\noindent
{\bf \gr 1.1. Theorem (Raynaud):}\rm\ {\sl The determinant $\det u$ of $u$ 
is not identically zero on $J^1$.} 

\pop {.5}
\par
In particular, one can consider the 
divisor $\theta:=\theta _{B}$ on $J^1$, which is the 
positive Cartier divisor locally generated by $\det u$. This is the {\it theta 
divisor} associated to the vector bundle $B$. By definition a point 
$a\in J^1(k)$ lies on the support of $\theta$ if and only if 
$H^0(X^1,B\otimes \Cal L_a)\neq 0$.

\pop {1}
\par
\noindent
{\bf \gr 2. $\bold p$-Rank of cyclic \'etale covers of degree prime to 
$\bold p$.}

\pop {.5}
\par
We use the same notation as in 1. We will only discuss in this section
the $p$-rank and the notion of {\it new-ordinariness} for cyclic 
covers of degree 
l := a prime integer distinct from $p$. This is the only case we 
use in this paper. For the general case of any integer prime to $p$ see 
[Ra-1], 2,  and [Ta], 3. 

\pop {.5}
\par
The absolute Frobenius morphism $F:X\to X$ induces
a semi-linear map $F:H^1(X,\Cal O_X)\to H^1(X,\Cal O_X)$, and we have a 
canonical decomposition: 
$$H^1(X,\Cal O_X)=H^1(X,\Cal O_X)^{\ss}\oplus 
H^1(X,\Cal O_X)^{\n}$$
where $H^1(X,\Cal O_X)^{\ss}$ is the semi-simple part 
on which $F$ is bijective, and $H^1(X,\Cal O_X)^{\n}$ is the nilpotent part on 
which $F$ is nilpotent. The $p$-{\it rank} $r_X:=r$ of $X$ is 
the dimension of the $k$-vector space $H^1(X,\Cal O_X)^{\ss}$. 
By duality, it is also the dimension of the subspace of 
$H^0(X,\Omega ^1_X)$ on which the Cartier operation $c$ is bijective
(cf. [Se], 10). The $p$-rank $r_X$ of $X$ is also the rank of the maximal 
pro-$p$-quotient $\pi _1^p(X)$ of the fundamental group $\pi _1(X)$ of $X$, 
which is known to be a finitely generated free pro-p-group (cf. [Sh]). If
$A$ is an abelian variety of dimension $d$ over $k$, then the rank of the \'etale part of the 
kernel of the morphism $[p]:A\to A$ of multiplication by $p$ is $p^h$,
where $0\le h \le d$ is the $p$-rank of $A$. 
The abelian variety $A$ is said to be 
{\it ordinary} if is has maximal $p$-rank equal to $d$, which is also 
equivalent to the fact that the Frobenius $F$ is bijective 
on $H^1(A,\Cal O_A)$. With the above notation, if $J=\Pic^0(X)$ is the 
jacobian variety of $X$, then it is well known that the $p$-rank of $X$ equals 
the $p$-rank of $J$.

\pop {.5}
\par
The relative Frobenius morphism $\pi :X\to X_1$ induces (because it is a radicial morphism) 
a ``canonical'' isomorphism $\pi _1(X)\to \pi _1(X_1)$ between fundamental groups (cf. [SGA-1], IX, 
Th\'eor\`eme 4. 10). In 
particular, for any {\bf prime} integer $l$, which is distinct from 
$p$, one has a one-to-one correspondence between $\mu_l$-torsors of 
$X$ and those of $X^1$. 
More precisely, the canonical homomorphism $H^1_{\et}(X^1,
\mu _l)\to H^1_{\et}(X,\mu _l)$, induced by $\pi $, is an isomorphism. 
Consider a $\mu _l$-torsor $f:Y\to X$ with $Y$ connected. By Kummer theory,
the torsor $f$ is given by an invertible sheaf $\Cal L$ of order $l$ on $X$,
and $Y:=\Spec (\oplus _{i=0}^{l-1}\Cal L^{\otimes i})$. There exists then an
invertible sheaf $\Cal L^1$ on $X^1$, of order $l$, such that if 
$f':Y^1\to X^1$ is the associated $\mu_l$-torsor, we have a cartesian diagram:

$$
\CD
Y     @>f>> X \\
@V\pi ' VV      @V\pi VV  \\
Y^1    @>f'>> X^1 \\
\endCD
$$

\pop {.5}
\par
Let $J_Y:=\Pic^0(Y)$ (resp. $J_X:=\Pic^0(X)$) denotes the Jacobian variety 
of $Y$ (resp. the Jacobian of $X$). The morphism $f:Y\to X$ induces a 
natural homomorphism $f^*:J_X\to J_Y$ between Jacobians, which has a 
finite kernel ($f^*$ is given by the pull-back of degree zero invertible 
sheaves). 
Let $J^{\new}:=J_{Y/X}$ denotes the quotient of $J_Y$
by the image $f^*(J_X)$ of $J_X$. The variety $J^{\new}$ is an abelian variety
of dimension $g_Y-g_X$, and $p$-rank equal to $r_Y-r_X$, it is called 
the {\it new part} of the Jacobian $J_Y$ of $Y$ with respect to the 
morphism $f$. 

\pop {.5}
\par
\noindent
{\bf \gr 2.1. Definition.}\ \rm The $\mu_ l$-torsor $f:Y\to X$ is said to be  
{\it new-ordinary}, if the new part $J^{\new}$ of the Jacobian of $Y$, with 
respect to the morphism $f$, is an ordinary abelian variety, i.e. if the 
equality $g_Y-g_X=r_Y-r_X$ holds.

\pop {.5}
\par
Raynaud's theory of theta divisors allows another important geometric 
interpretation of new-ordinariness, which we explain below. This interpretation
has allowed significant recent progress in the study of fundamental groups 
of curves in positive characteristics.

\pop {.5}
\par
There exists an isomorphism 
$H^1(J_Y,\Cal O_{J_Y})\simeq  H^1(Y,\Cal O_Y)$ (cf. [Se-1], 
VII, th\'eor\`eme 9), and $H^1(Y,\Cal O_Y)=H^1(X,f^*\Cal O_Y)=
H^1(X,\oplus_{i=0}^{l-1}\Cal L^{\otimes i})$, from which we deduce that 
$H^1(J^{\new},
\Cal O_{J^{\new}})\simeq H^1(X,\oplus_{i=1}^{l-1}
(\Cal L)^{\otimes i})$. Note that the above identifications are compatible 
with the action of Frobenius. Hence, the kernel of Frobenius on $H^1(J^{\new},
\Cal O_{J^{\new}})$ is isomorphic to the kernel of Frobenius acting on
$H^1(X,\oplus_{i=1}^{l-1}\Cal L^{\otimes i})$.
On the other hand, as $f'$ is \'etale, we have $(f')^*(B_X)=B_Y$. Thus, also 
$(f')_*(B_Y)=B_X\otimes (f')_*(\Cal O_{Y^1})=\oplus_
{i=0}^{l-1}(B_X\otimes (\Cal L^1)^{\otimes i})$. Now, by duality, the 
kernel of the Frobenius acting on
$H^1(X^1,\oplus_{i=1}^{l-1}{\Cal L^1}^{\otimes i})$ is 
isomorphic to the kernel of the Cartier operator acting on 
$H^0(X^1,\pi _*\Omega ^1_X \otimes (\oplus_{i=1}^{l-1}
(\Cal L^1)^{\otimes i}))$, which is $\oplus_{i=1}^{l-1}
H^0(X^1,B_X\otimes (\Cal L^1)^{\otimes i})$. Thus, we see that 
the above $\mu_ l$-torsor $f:Y\to X$ is new ordinary, if and only if 
the Frobenius $F$ is injective (hence bijective) on $H^1(J^{\new},
\Cal O_{J^{\new}})$, i.e. if and only if
$H^0(X^1,B\otimes {({\Cal L}^1)}^{\otimes i})=0$ for all $i\in \{1,...,l-1\}$. 
Finally, this last statement is equivalent, by the very definition of the 
theta divisor $\theta _X$ associated to the vector bundle $B_X$, 
to the following:

\pop {.5}
\par
\noindent
{\bf \gr 2.2. Proposition:}\ \rm {\sl The $\mu_ l$-torsor $f:Y\to X$ is   
new-ordinary, if and only if the subgroup $<{\Cal L}^1>$, 
generated by ${\Cal L}^1$ in $J^1$, intersects 
the support of the theta divisor $\theta _X$ at most 
at the zero point $0_{J^1}$ of $J^1$}.

\pop {.5}
\par
Using the above interpretation of new-ordinariness, and with an input form 
intersection theory, one can prove that for $l>>0$ ``most'' 
$\mu_l$-torsors are new-ordinary. More precisely, one has the following 
result which is essencially due to Serre and Raynaud 
(see [Ra], th\'eor\`eme 4.3.1, and [Ta], corollary 3.10, for a proof):

\pop {.5}
\par
\noindent
{\bf \gr 2.3. Theorem:}\rm\ {\sl There exists a constant $c$, depending 
only on $g$ and $p$, such that for each prime integer $l\neq p$, 
the set of elements of $J[l](k)$ whose corresponding 
$\mu_l$-torsor is not new-ordinary, has cardinality 
$\le c(l-1)l^{2g-2}$ (Here $J[l]$ denotes the kernel of multiplication by 
$l$ in $J$). Moreover, one can take 
$c=(p-1)3^{g-1}g!$}.

\pop {.5}
\par
In particular, if $l>>0$ we can find an element
of $J[l](k)$ such that the corresponding $\mu_l$-torsor is new-ordinary, since $\Card J[l](k)=l^{2g}$.

\pop {1}
\par
\noindent
{\bf \gr 3. Equi-characteristic deformation of generalized Prym varieties.}

\pop {.5}
\par
In this section, we state the theorem of Tamagawa, on the local infinitesimal 
Torelli problem for generalized Prym varieties. This Theorem is an essential 
tool in the proof of the main result of this paper.

\pop {.5}
\par
Let $k$ be an algebraically closed field, of arbitrary characteristic. 
Denote by $\Cal C_k$ the category
of artinian local rings, with residue field $k$. For a proper and smooth 
$k$-variety $X_0$, one defines the (equi-characteristic) deformation 
functor $M_{X_0}$ of $X_{0}$, to be the functor: 
$$M_{X_0}:\Cal C_k\to (\Sets)$$ 
which to an element $R$ of $\Cal C_k$, associates the set of isomorphism 
classes of pairs $(X,\varphi)$, where $X$ is a proper and smooth $R$-scheme, 
and $\varphi$ is an isomorphism $X\times _Rk\simeq X_0$. The functor $M_{X_0}$
is well understood in the case where $X_0$ has dimension $1$, and genus $g> 1$,
(resp. if $X_0$ is an abelian variety of dimension $d$). In this case 
the functor $M_{X_0}$ is pro-representable by a ring of formal power series
of $3g-3$, (resp. $d^2$), variables over $k$. We will be mainly interested 
in this two cases.

\pop {.5}
\par
Now, assume that $X_0$ is a proper, connected, and smooth algebraic 
curve over $k$, with genus $g\ge 2$. Let $l$ be a prime integer 
distinct from the characterisic of $k$, and let 
$f_0:Y_0\to X_0$ be a $\mu_l$-torsor, with 
$Y_0$ connected. The torsor $f_0$ corresponds to an element 
${\Cal L}_0\in J_0[l](k)$, where
$J_0[l](k)$ denotes the $k$-subgroup, of $l$-torsion points, in the jacobian 
$J_0$ of $X_0$. Let $J_0^{\new}$ be the new part of the jacobian of $Y_0$, 
with respect to the morphism $f_0$. Then, for any element $R$ of $\Cal C_k$, 
there is a natural map: 
$$T_{{\Cal L}_0}(R):M_{X_0} (R)\to M_{J_0^{\new}}(R)$$ 
defined as follows: 
let $(X,\varphi)$ be an element of $M_{X_0}(R)$. The $\mu_l$-torsor
$f_0:Y_0\to X_0$ lifts uniquely, by the theorems of lifting of \'etale covers
(cf. [SGA-1], I, 8), to a $\mu_l$-torsor $f:Y\to X$. 
Let $J_X:=\Pic^0(X)$ (resp. 
$J_Y:=\Pic^0(Y)$) be the relative jacobian of $X$ (resp. the relative jacobian 
of $Y$), which is an 
abelian scheme over $R$, and let
$f^*:J_X\to J_Y$ be the natural homomorphism, which is 
induced by the pull back of invertible sheaves. 
Define $J^{\new}:=J_Y/f^*(J_X)$, to be the quotient of
$J_Y$ by the image $f^*(J_X)$ of $J_X$. Then $J^{\new}$ is an abelian 
$R$-scheme, and there exists a natural isomorphism 
$\psi:J^{\new}\times _R k\simeq J_0^{\new}$, which is induced by $\varphi$. 
Thus, the pair
$(J^{\new},\psi)$ is an element of $M_{J_0^{\new}}(R)$, which we define 
to be the image under $T_{\Cal L_0}(R)$ of $(X,\varphi)$ 
in $M_{J_0^{\new}}(R)$. The infinitesimal 
Torelli problem asks whether or not the above natural map: 
$$T_{{\Cal L}_0} : M_{X_0}\to M_{J_0^{\new}}$$
is an immersion. More precisely, let $k[\epsilon]$ ($\epsilon ^2=0$) be the 
ring of dual numbers on $k$. Then, the question is whether the natural map:
$T_{{\Cal L}_0}(k[\epsilon]) : M_{X_0} (k[\epsilon])\to 
M_{J_0^{\new}}(k[\epsilon])$, between tangent spaces, is injective. 
This is also equivalent asking, whether if $A$ (resp. $B$) is the 
pro-representing object of the functor $M_{X_0}$ (resp. of 
the functor $M_{J_0^{\new}}$), then the natural homomorphism $B\to A$,
induced by $T_{\Cal L_0}$, is surjective. If this is the case, it would 
in particular imply the following: for every element $R\in \Cal C_k$, 
and $(X,\varphi)\in M_{X_0}(R)$, if the image $(J^{\new},\psi)$ of
$(X,\varphi)$ in $M_{J_0^{\new}}(R)$, via the map $T_{{\Cal L}_0}(R)$,
is a trivial deformation of $J_0^{\new}$, then $X$ is a trivial 
deformation of $X_0$. This is 
a generalization of the classical infinitesimal Torelli problem, which asks 
whether the natural map from the space of deformations of the curve $X_0$, 
to the space of deformations of its jacobian, is an immersion, in which case
one knows that the answer is yes, if the curve $X_0$ is not hyperelliptic. 
Concerning the above generalization, Tamagawa proves the following:

\pop {.5}
\par
\noindent
{\bf \gr 3.1. Theorem (Tamagawa):}\rm
\ {\sl Let $X_0$ be a proper, connected, and 
smooth algebraic curve, over an algebraically closed field $k$, with genus 
$g\ge 2$. Let $d_{X_0}:=\{\min \deg (f)\ /\ f:X_0\to \Bbb P^1_k$ 
non constant$\}$,
be the gonality of $X_0$ (cf. [Ta-1], 1). Assume that $d_{X_0}\ge 5$.
Then, there exists a constant $c_1$, which depends only on $X_0$, and 
such that for each prime integer $l\neq\char(k)$, 
the subset of the elements $\Cal L_0$ of $J[l](k)$, 
such that the corresponding natural map $T_{{\Cal L}_0} : M_{X_0}\to 
M_{J_0^{\new}}$ is not an immersion, has cardinality $\le c_1l^{2g-2}$.} 

\pop {.5}
\par
In particular, if $l>>0$, then one can find an element $\Cal L_0\in J[l](k)$,
such that the corresponding map: $T_{{\Cal L}_0} : M_{X_0}\to M_{J_0^{\new}}$ 
is an immersion. For a proof of the above result see [Ta-1], corollary 4.16. 
Concerning the gonality of curves, Tamagawa also proves the following:

\pop {.5}
\par
\noindent
{\bf \gr 3.2. Theorem (Tamagawa):}\rm\ {\sl Let $X$ be a 
proper, connected, and smooth algebraic curve, over an algebraically closed 
field $k$, with genus $g\ge 2$. Then, there 
exists an \'etale cover $f:Y\to X$, such that the gonality $d_Y$ of $Y$ 
satisfies $d_Y\ge 5$. Moreover, the cover $f$ can be chosen to be a 
composition of two cyclic \'etale covers, of (suitable) degree prime to the 
characteristic of $k$.}

\pop {.5}
\par
For the proof of Theorem 3.2, combine Theorem 2.7, Proposition 2.14, 
and Corollary 2.19 from [Ta-1].

\pop {.5}
\par 
Combining both the Theorems 3.2 and 2.3 above, we obtain the following
result, which we will use in the proof of our main Theorem in section 
6. This result was also used by Tamagawa, in [Ta-1], in order 
to prove Theorem 6.1 in that paper.

\pop {.5}
\par
\noindent
{\bf \gr 3.3. Theorem:}\rm\ {\sl Let $X$ be a 
proper, connected, and smooth algebraic curve, of genus  $g\ge 2$, over an 
algebraically closed field $k$ of characteristic $p>0$. Assume that the 
gonality $d_X$ of $X$ is $\ge 5$. Then, if $l\neq p$ is a prime integer,
such that $l>1+c_1+(p-1)3^{g-1}g!$, where $c_1$ is 
the constant in theorem 3.2, there exists
a non zero element $\Cal L\in J[l](k)$, 
such that the following two conditions are satisfied:
\par
(i)\ The $\mu_l$-torsor $f:Y\to X$ corresponding to $\Cal L$ is new ordinary.
\par
(ii)\ The natural map $T_{\Cal L} : M_{X}\to M_{J^{\new}}$ is an immersion.}

\pop {1}
\par
\noindent
{\bf \gr 4. Complete families of ordinary abelian varieties in characteristic $p>0$.}

\pop {.5}
\par
In this section, we state the theorem of Raynaud, Szpiro, and Moret-Bailly, 
on the isotriviality of complete families of ordinary abelian varieties, 
in characteristic $p>0$. This result will be used in section 6. In all what follows we 
fix a prime integer $p>0$.

\pop {.5}
\par
Let $g\ge 1$, and $d\ge 1$ be integers. Let $\Cal A_{g,d}\to \Spec\Bbb F_p$ denotes the 
coarse moduli scheme of polarized abelian varieties (with a degree $d$ polarization), of dimension $g$, 
in characteristic $p$. The scheme  $\Cal A_{g,d}$ is a quasi-projective variety of dimension $g(g+1)/2$, and 
has the following property: for every scheme $S$ in characteristic $p$, 
and $X\to S$ a polarized abelian $S$-scheme (with a degree $d$ polarization), of relative 
dimension $g$, there exists a unique natural map $S\to {\Cal A}_{g,d}$, defined by 
the family $X\to S$. This map sends a point $s\in S$, to the 
moduli point corresponding to the fibre $X_s\to \Spec k(s)$, of $X$ above 
the point $s$.

\pop {.5}
\par
\noindent
{\bf \gr 4.1. Definition:}\rm\ Let $k$ be a field of characteristic $p>0$.
Let $S$ be a normal and integral $k$-variety, with generic point $\eta$, and let $X\to S$ be an abelian $S$-scheme of relative 
dimension $g$. Assume that the generic fibre $X_{\eta}\to \Spec k(\eta)$ of $X$, has a polarization of degree $d$. 
Then the family $X\to S$ is called {\it isotrivial}, if the
natural map $\Spec k(\eta)\to {\Cal A}_{g,d}$, defined by $X_{\eta}\to \Spec k(\eta)$, factorizes through
$\Spec k(\eta)\to \Spec k'\to\Spec {\Cal A}_{g,d}$, where $k'$ is a finite extension of $k$.

\pop {.5}
\par
\noindent
{\bf \gr 4.2. Theorem (Raynaud, Szpiro, Moret-Bailly):}\rm\ {\sl  Let $k$ be a field of characteristic $p>0$. 
Let $S$ be a normal and integral $k$-variety, which is {\bf projective}. Let $X\to S$ be an abelian $S$-scheme, 
of relative dimension $g$. Assume that all the geometric fibres $X_{\bar s}\to \Spec k({\bar s})$, 
of $X$ over $S$, are ordinary abelian varieties. Then the family $X\to S$ is isotrivial.}

\pop {.5}
\par
Next, we would like to explain the stratification, in characteristic $p$, of 
the moduli space of polarized abelian varieties, by the $p$-rank. Let $k$ be an 
algebraically closed field of characteristic $p$. For each fixed integer $0\le f\le g$,
let $\Cal V_f\subset \Cal A_g\times _{\Bbb F_p}k$ be the subset of points 
corresponding to abelian varieties having a $p$-rank $\le f$. Concerning the 
subsets $\Cal V_f$, we have the following:

\pop {.5}
\par
\noindent
{\bf \gr 4.3. Theorem:}\rm\ {\sl For each fixed integer $0\le f\le g$, and any 
algebraically closed field $k$ of characteristic $p$, let $\Cal V_f\subset \Cal A_g\times _{\Bbb F_p}k$
be the subset of points corresponding to abelian varieties 
having a $p$-rank $\le f$. Then, the subset $\Cal V_f$ is a closed subscheme of
$\Cal A_g\times _{\Bbb F_p}k$, and every irreducible component of 
$\Cal V_f$ has dimension equal to $g(g+1)/2-g+f$. Moreover, the closed 
subscheme $\Cal V_0$ is a complete subvariety of 
$\Cal A_g\times _{\Bbb F_p}k$, 
of dimension $g(g-1)/2$}.
\pop {.5}
\par
For the proof of the fact that $\Cal V_f$ is closed see [Oo-1], Corollary 1.5. 
For the statement concerning the dimension of the above strata 
see [No-Oo], Theorem 4.1. Finally, for the fact that $\Cal V_0$ is complete, 
see [Oo-1], proof of Theorem 1.1 a).

\pop {1}
\par
\noindent
{\bf \gr 5. Complete families of curves with constant $p$-rank in characteristic $p>0$.}
\pop {.5}
\par
In this section, we want to extend the argument of Oort, in [Oo-1], 
in order to show the existence,
for every integer $g>2$, of a non-isotrivial complete family of smooth and 
proper curves of genus $g$, with constant $p$-rank equal to $0$. 
Before stating the main result, we will explain, and adopt some notation. 
In all what follows we fix a prime integer $p>0$.

\pop {.5}
\par
Let $g\ge 2$ be an integer, and let $\Cal M_g\to \Spec \Bbb F_p$ be 
the coarse moduli scheme, of projective, smooth, and irreducible curves of 
genus $g$, in characteristic $p$. The
scheme $\Cal M_g$ is a quasi-projective irreducible 
variety, of dimension $3g-3$. 
Let $\overline {\Cal M}_g\to \Spec \Bbb F_p$ be the Deligne-Mumford 
compactification of $\Cal {M}_g$,
which is a coarse moduli scheme, of projective, and stable curves of genus 
$g$, in characteristic $p$. The scheme $\overline {\Cal M}_g$ is 
an irreducible projective variety, which contains 
$\Cal M_g$ as an open subscheme (cf. [De-Mu] for more details. 
We will also consider $\Cal {M}'_g\to \Spec \Bbb F_p$,
which is the coarse moduli scheme of projective, and stable curves of genus 
$g$, in characteristic $p$, whose jacobian is an abelian variety 
(these are exactly the  projective, and stable curves of genus $g$, 
whose configuration is a tree like cf. [Bo-Lu-Ra], 9, corollary 12).

\pop {.5}
\par
Let $S$ be a scheme of characteristic $p$. By a smooth relative 
(or a family of) curve(s) $f:X\to S$ over $S$, of genus $g$, we mean 
that $f$ is a proper and smooth equidimensional morphism, 
with relative dimension $1$, whose fibres are 
curves of genus $g$. We say that the relative curve 
$f:X\to S$ is {\it complete}, if $S$ is a complete, and irreducible 
variety, in characteristic $p$. The moduli scheme ${\Cal {M}}_g$ has the following property: if 
$f:X\to S$ is a smooth relative curve over $S$, of genus $g$, then there is a natural map
$S\to{\Cal {M}}_g$, which is uniquely determined by $f$. This map sends a point $s\in S$, to the 
moduli point corresponding to the fibre $X_s\to \Spec k(s)$, of $X$ above the point $s$. We also have the 
following:

\pop {.5}
\par
\noindent
{\bf \gr 5.1. Proposition:}\rm\ {\sl Let $k$ be a field of characteristic $p$, and let $\Cal S$ be a subvariety 
of $\Cal M_g\times_{\Bbb F_p}k$. Then, there exists a finite \'etale cover 
$h:\Cal S'\to \Cal S$, and a smooth 
relative $\Cal S'$-curve $f':X'\to \Cal S'$ of genus $g$, such that 
the natural morphism $\Cal S'\to \Cal {M}_g$, induced by $f'$, 
factorizes $\Cal S'\to S\to \Cal {M}_g$ through $h$.}

\pop {.5}
\par
\noindent
{\bf \gr Proof:}\rm\ Standard, by passing to the fine moduli scheme of smooth,
projective, and irreducible curves of genus $g$, with a symplectic level $n$-structure.

\pop {.5}
\par
\noindent
{\bf \gr 5.2. Definition:}\rm\ Let $S$ be a scheme of characteristic $p$,
and let $f:X\to S$ be a smooth relative $S$-curve, of genus $g$. The curve
$f:X\to S$ is said to be {\it isotrivial}, if the corresponding map 
$S\to {\Cal {M}}_g$, has an image which consists of a point.

\pop {.5}
\par
\noindent
{\bf \gr 5.3. Definition:}\rm\ Let $S$ be a scheme of characteristic $p$,
and let $f:X\to S$ be a smooth relative $S$-curve of genus $g$. 
Let $0\le r\le g$ be an integer. We say that the family $f:X\to S$ has 
constant $p$-rank, equal to $r$, if each geometric fibre 
$X_{\bar s}\to \Spec k({\bar s})$ of $f$ has a $p$-rank equal to $r$.

\pop {.5}
\par
In [Oo-1], Oort showed in the proof of theorem 1.1. b), the existence over any 
algebraically closed field $k$ of characteristic $p$, of a complete curve 
contained in $\Cal M_3\times _{\Bbb F_p}k$, which is contained in 
the locus of curves having $p$-rank equal to $0$. This indeed corresponds 
by Proposition 5.1, to a non-isotrivial complete family of smooth curves, 
of genus $3$, with constant $p$-rank equal to $0$. 
Next, we want to extend Oort's argument, in order to prove the following:

\pop {.5}
\par
\noindent
{\bf \gr 5.4. Theorem:}\rm\ {\sl Let $g\ge 3$ be an integer. Let
$k$ be an algebraically closed field, of characteristic $p$. Then 
$\Cal M_g\times _{\Bbb F_p}k$ contains a complete irreducible curve, which is 
contained in the locus of curves having $p$-rank equal to $0$.}

\pop {.5}
\par
\noindent
{\bf \gr Proof:}\rm\ The proof consists in considering the Torelli map, 
together with a dimension argument. More precisely, let: 
$$t:\Cal M'_g\times _{\Bbb F_p}k\to \Cal A_{g,1}\times _{\Bbb F_p}k$$
be the Torelli morphism, which sends the class of a stable curve, whose 
jacobian is an abelian variety, to the class of its jacobian, endowed
with its canonical principal polarization coming from the theta divisor.
Torelli's theorm (cf. [Mi], 12) says that the map $t$ is injective on 
geometric points. In particular, the image $\Cal J'_g:=t(\Cal M'_g\times _{\Bbb F_p}k)$ (resp.
$\Cal J_g:=t(\Cal M_g\times _{\Bbb F_p}k)$) of 
$\Cal M'_g\times _{\Bbb F_p}k$ (resp, of $\Cal M_g\times _{\Bbb F_p}k$), which is called the 
{\it jacobian locus} (resp. {\it open jacobian locus}), 
is a subvariety of $\Cal A_{g,1}$ of dimension $3g-3$. 
Moreover, $\Cal J'_g$ is closed in $\Cal A_{g,1}$ (cf. [Mu], lecture IV, p. 74).
For each fixed integer $0\le f\le g$, let $\Cal V_{f}$ be 
the closed subscheme of $\Cal A_{g,1}\times _{\Bbb F_p}k$ as defined in 4.3. 
Then, every irreducible 
component of $\Cal J'_{f,g}:=\Cal J'_g\cap \Cal V_{f}$, has dimension 
at least $3g-3-g+f=2g-3+f$ (cf. [Oo-1], lemma 1.6). Moreover, 
it is well known, that every irreducible 
component of $\Cal J'_{0,g}:=\Cal J'_g\cap \Cal V_{0}$ has dimension
$2g-3$ (cf. [Fa-Lo], 11, p. 36). We claim, that 
$\{\Cal J'_g-\Cal J_g\}\cap \Cal V_0$, has codimension at least
two, in  $\Cal J'_{0,g}:=\Cal J'_g\cap \Cal V_{0}$. 
Indeed, $\{\Cal J'_g-\Cal J_g\}\cap \Cal V_0$
is contained in the images of $\Cal J_{0,g_1}\times...\times
\Cal J_{0,g_t}$, for all possible family of integers $\{g_1,...,g_t\}$, 
such that $g_1+...+g_t=g$, with $t\ge 2$, via the natural morphism 
$\Cal J'_{0,g_1}\times...\times\Cal J'_{0,g_t}   \to \Cal V_{0,g}$. 
Now, counting the dimension of $\Cal J'_{0,g_1}\times...\times\Cal J'_{0,g_t}$,
which is $2g_1-3+...+2g_t-3 < 2g-5$, 
we conclude that $\{\Cal J'_g-\Cal J_g\}\cap \Cal V_0$,
has codimension at least two in $\Cal J'_{0,g}:=\Cal J'_g\cap \Cal V_{0}$.
\par
Finally, since $\Cal V_0$ is projective, we can find a closed immersion 
$\Cal V_0\to \Bbb P^N_k$ into a projective $k$-space, of suitable dimension.
Further, we can find a general linear subspace $L$ of $\Bbb P^N_k$, 
of suitable dimension, such that $L\cap \Cal J'_{0,g}$ is a 
(necessarily complete) curve $S'$,
and such that  $L\cap \{\Cal J'_{0,g}-\{\Cal J_{g}\cap \Cal V_0\}\}$ is empty 
(cf. [Da-Sh], II, Chapter 3, 1.2). Now the inverse image $S:=t^{-1}(S')$, 
of $S'$ via the Torelli map $t$, is a complete curve, contained in 
$\Cal M_g \times_{\Bbb F_p} k$, and which by construction, is contained in the 
locus of curves having $p$-rank equal to $0$.

\pop {.5}
\par
\noindent
{\bf \gr 5.5. Corollary:}\rm\ {\sl Let $g\ge 3$ be an integer. Let
$k$ be an algebraically closed field, of characteristic $p$. Then there
exists a complete, and smooth algebraic $k$-curve $S$, and a non iso-trivial
smooth $S$-curve $f:X\to S$, of genus $g$, 
with constant $p$-rank equal to $0$.}

\pop {.5}
\par
\noindent
{\bf \gr 5.6. Remark:}\rm\ It is tempting to try to construct an $S$-curve $X$,
as in 5.5, for special values of $g$, by considering a Galois cover 
$f:X\to Y$, with group
$\Bbb Z/p\Bbb Z$, where $Y$ is a ruled surface over $S$, and such that $f$ is
\'etale outside an $S$-section of $Y$. The Deuring-Shafarevich formula,
comparing the $p$-rank in Galois $p$-covers, would then imply 
that all fibres of $X$ over $S$ have $p$-rank equal to $0$. 
However, by a result of Pries, all such covers $f:X\to Y$ are 
necessarily isotrivial (cf. [Pr], Theorem 3.3.4).

\pop {.5}
\par
\noindent
{\bf \gr 5.7. Question:}\rm\  Let $g\ge 3$ be an integer. Let
$k$ be an algebraically closed field, of characteristic $p$. Let $\tilde r$
be the maximum of the integers $r$, such that $\Cal M_g\times _{\Bbb F_p}k$ 
contains a complete curve $S$, which is contained in $\Cal V_r-\Cal V_{r-1}$. 
We have $0\le \tilde r<g$. What
is the value of $\tilde r$? Does the value of $\tilde r$ depends only on $g$?
Or does it depend on $p$ as well?

\pop {1}
\par
\noindent
{\bf \gr 6. Complete families of curves with a given fundamental group in characteristic $p>0$.}

\pop {.5}
\par
This is our main section, in which we prove the main result of this 
paper, which asserts that complete families of curves with a 
constant geometric fundamental group are 
isotrivial. In all what follows we fix a prime intger $p>0$.

\pop {.5}
\par
First, we will explain how to define the fundamental group of points, in 
the moduli space of curves. Let $\Cal M_g\to \Spec \Bbb F_p$ be the coarse 
moduli scheme of smooth, and projective curves of genus $g$, in 
characteristic $p>0$. 
Let $k$ be an algebraically closed field of characteristic $p>0$. 
For a geometric point $\overline x\in \Cal M_g(k)$, let 
$C_{\overline x}\to \Spec k$ be a smooth curve, of genus $g$, which is 
classified by $\overline x$, and let
$x\in \Cal M_g$ be the point such that $\overline x:\Spec k\to \Cal M_g$ 
factors through $x$. We
define the {\it geometric fundamental group} 
$\pi _1(x):=\pi _1(C_{\overline x})$,
of the point $x$, as the \'etale fundamental group of the
curve $C_{\overline x}$ (we assume of course the choice of a base point).
We remark that the structure of $\pi _1(x)$, as a profinite group, depends 
only on the point $x$, and not on the concrete geometric point 
$\overline x\in \Cal M_g(k)$ used to define it. 
Indeed, first, if $\overline {k(x)}$ is an algebraic closure 
of the residue field $k(x)$ at $x$, and $C_x$ is the curve classified by
$\Spec \overline {k(x)}\to \Cal M_g$, then 
$C_{\overline x}\simeq C_x\times _{\overline {k(x)}} k$, is the base change of 
$C_x$ to $k$. Hence, $\pi _1(C_{\overline x})\simeq \pi _1(C_x)$,
by the geometric invariance of the fundamental group for proper varieties
(cf. [SGA-1], X, Corollaire 1.8). Second, the isomorphy type of 
$C_x$ as an $\Bbb F_p$-scheme 
does not depend on the choice of $\overline {k(x)}$: a geometric point of 
$M_g$ dominating $x$. We thus have a map:
$$\pi _{1,\geom}:\Cal M_g\to \{\ProfGrp\}$$
$$x\to \pi _1(x)$$

\pop {.5}
\par
Next, we recall the specialization theory of 
Grothendieck for fundamental groups. Let $y\in \Cal M_g$ be a point which 
specializes into the point $x\in \Cal M_g$. Then Grothendieck's 
specialization theorem shows the existence, of
a surjective continuous homomorphism $\Sp:\pi _1(y)\to \pi_1(x)$
(cf. [SGA-1], X). In particular, if $\eta$ is the generic point of $\Cal M_g$, 
then $C_{\eta}$ is the generic curve 
of genus $g$, and every point $x$ of $\Cal M_g$ is a specialization of $\eta$.
Hence, for every $x\in \Cal M_g$, there exists a surjective homomorphism 
$\Sp_x:\pi _1(\eta)\to \pi _1(x)$. For every such an $x$, we fix such a 
map once for all. In particular, if $y$ specializes to $x$, 
we also fix a surjective homomorphism 
$\Sp_{y,x}:\pi _1(y)\to \pi_1(x)$, 
such that $\Sp _{y,x}\circ \Sp _y=\Sp _x$.

\pop {.5}
\par
\noindent
{\bf \gr 6.1. Definition.}\rm\ Let $\Cal S\subset \Cal M_g$ be a subscheme of 
$\Cal M_g$. We say that the (geometric) fundamental group $\pi_1$ is
{\it constant} on $\Cal S$, if for any two points $x$ and $y$ of $\Cal S$, such
that $y$ specializes in $x$, the corresponding specialization 
homomorphism $\Sp_{y,x}:\pi _1(y)\to \pi_1(x)$ is an isomorphism. We say
that $\pi_1$ is {\it not constant} on $\Cal S$, if the contrary holds, 
namely: there exists two points $x$ and $y$ of $\Cal S$, such
that $y$ specializes in $x$, and such that the corresponding specialization 
homomorphism $\Sp_{y,x}:\pi _1(y)\to \pi_1(x)$ is not an isomorphism.

\pop {.5}
\par
For every field $k$ of characteristic $p$, we define in a similar way, the
geometric fundamental group of points in $\Cal M_g\times _{\Bbb F_p}k$, as well as
the notion of a subvariety $\Cal S \subset \Cal M_g\times _{\Bbb F_p}k$, 
on which the geometric fundamental group $\pi_1$ is constant.

\pop {.5}
\par
\noindent
{\bf \gr 6.2. Definition.}\rm\ Let $S$ be a scheme of characteristic $p$, and
let $f:X\to S$ be a relative smooth $S$-curve of genus $g$. We say 
that the (geometric) fundamental group $\pi_1$, is
{\it constant} on the family $f$, if for any two points $\eta$ and $s$ of 
$S$, such that $\eta$ specializes in $s$, the corresponding specialization 
homomorphism $\Sp_{\eta,s}:\pi _1(X_{\bar {\eta}})\to \pi_1(X_{\bar s})$ 
is an isomorphism, where $X_{\bar {\eta}}:=X\times _S{k(\bar {\eta})}$
(resp. $X_{\bar s}:=X\times _S{k(\bar s)}$) is the geometric fibre of 
$X$ over the point $\eta$, (resp. the geometric fibre of 
$X$ over the point $s$). If the above condition doesn't hold, we say
that the fundamental group $\pi_1$ is not constant on the family $f$.

\pop {.5}
\par
Let $S$ be a scheme of characteristic $p$, and let $f:X\to S$ be a relative 
smooth $S$-curve of genus $g$. It is clear, that if the fundamental group
is constant on the family $f$, and if $h:S\to \Cal M_g$ is the map induced 
by the family $f$, then the fundamental group is also constant, 
on the image $h(S)$ of $S$ in $\Cal M_g$. It is quite natural to ask the following question:

\pop {.5}
\par
\noindent
{\bf \gr 6.3. Question.}\rm\ Let $k$ be an algebraically closed field of 
characteristic $p$. Does $\Cal M_g\times _{\Bbb F_p}k$ contain, $k$-subvarieties of positive dimension $>0$, 
on which $\pi_1$ is constant?

\pop {.5}
\par
Our main result is the following:

\pop {.5}
\par
\noindent
{\bf \gr 6.4. Theorem (Main Result).}\rm\ {\sl Let  $k$ 
be an algebraically closed field, of characteristic $p$. 
Let $\Cal S\subset \Cal M_g\times _{\Bbb F_p}k$ 
be a {\bf complete} $k$-subvariety of $\Cal M_g\times _{\Bbb F_p}k$. Then, the 
fundamental group $\pi_1$ is {\it not constant} on $\Cal S$.

\pop {.5}
\par
For the proof of Theorem 6.4, it is clear, that one can reduce 
to the case where $S$ is a complete, and irreducible curve. 
The proof of Theorem 6.4, then follows easily, by using 5.1, from the 
following Theorem 6.5:

\pop {.5}
\par
\noindent
{\bf \gr 6.5. Theorem.}\rm\ {\sl Let $k$ be an algebraically closed field, of 
characteristic $p$. Let $S$ be a smooth, complete, and irreducible 
$k$-curve. Let $f:X\to S$ be a {\bf non-isotrivial} smooth family of curves 
of genus $g\ge 2$. Then, the fundamental group $\pi_1$ is not constant on 
the family $f$. In particular, if $h:S\to \Cal M_g\times _{\Bbb F_p}k$ 
is the map defined by $f$, then the fundamental group $\pi_1$ is not 
constant on the image $h(S)$ of $S$ in $\Cal M_g\times _{\Bbb F_p}k$}.

\pop {.5}
\par
In the process of proving Theorem 6.5 we prove, in fact, the following 
more precise result:

\pop {.5}
\par
\noindent
{\bf \gr 6.6. Theorem.}\rm\ {\sl Let $k$ be an algebraically closed field, of 
characteristic $p$. Let $S$ be a smooth, complete, and irreducible $k$-curve.
Let $f:X\to S$ be a {\bf non-isotrivial} smooth family of curves of 
genus $g\ge 2$. Then, there exist a finite \'etale cover 
$S'\to S$, an \'etale cover $Y'\to X':=X\times _SS'$ of degree prime to $p$, and
a closed point $s_0\in S'$ such that the $p$-rank of the geometric fibre $Y'_{k(\bar s_0)}\to k(\bar s_0)$, 
of $Y'$ above the point $s_0$, is strictly smaller, that the $p$-rank
of the generic geometric fibre  $Y'_{k(\bar \eta)}\to k(\bar \eta)$, of $Y'$
above the generic point $\eta$ of $S'$.}

\pop {.5}
\par
First, we start with the following lemmas:

\pop {.5}
\par 
\noindent
{\bf \gr 6.7. Lemma/Definition.}\rm\ {\sl Let $k$ be an algebraically 
closed field, of characteristic $p$. Let $S$ be a smooth, complete, and 
irreducible $k$-curve. Let $f:X\to S$ be a smooth family of 
curves, of genus $g\ge 2$. Let $s$ be a closed point of $S$, and let 
$f_s:Y_s\to X_s:=X\times _Sk(s)$ be a $\mu_n$-torsor, over the fibre of $X$ 
above the point $s$, where $n$ is coprime to $p$. 
Then, there exists a positive integer $d$, such that if $n$ is coprime to $d$, then there exists
a finite \'etale cover $h:S'\to S$, a $\mu_n$-torsor 
$f':Y'\to X':=X\times _S S'$, and a closed point $s'\in S'$, with $h(s')=s$, such that the fibre 
$f'_{s'}:Y'_{s'}:=Y'\times _{S'} k(s')\to X'_{s'}
:=X'\times _{S'} k(s')$, of the torsor $f'$, above the point $s'\in S'$, 
coincides with the given cover $f_s:Y_s\to X_s$. We call such a pair $(f',h)$ 
a {\bf good lifting} of the cover $f_s:Y_s\to X_s$.}

\pop {.5}
\par
\noindent
{\bf \gr Proof.}\rm\  In the case where the morphism $f$ has a section, the above 
lemma follows easily from the homotopy exact sequence of fundamental groups in
[SGA-1], XIII, Proposition 4.3 (see Lemma 4.3.1 in Loc.cit). In the general case, let $x$ be a 
closed point of the genric fibre $X_{\eta}$ of $X$ over $S$, and let $Z$ be the schematic closure 
of $x$ in $X$. Denote by $Y$ be the normalization of $Z$. The canonical morphism $Y\to S$ 
is finite of degree $d$, it is a
``multisection'' of $f$ of degree $d$. Assume, further, that the integer $n$ is coprime to $d$. 
The sheaf $R^1f_{*}(\Bbb Z /n\Bbb Z)$
is locally constant on $S_{\et}$. In particular, there exists a finite \'etale cover $h:S'\to S$,
such that $R^1f_{*}(\Bbb Z /n\Bbb Z)/S'$ is constant. We denote by
$Y'\to X':=X\times_SS'\to S'$, a multisection of $f':X'\to S'$ above $Y$.
The Leray spectral sequence in \'etale 
cohomology, with respect to the morphism $f':X':=X\times _SS'\to S'$, and the constant sheaf
$\Bbb Z/n\Bbb Z$, gives rise to an exact sequence of terms of low degree: $0\to H^1(S',f'_*
(\Bbb Z/n\Bbb Z))\to H^1(X',\Bbb Z/n\Bbb Z)\to H^0(S',R^1f'_*
(\Bbb Z/n\Bbb Z))\to H^2(S',f'_*(\Bbb Z/n\Bbb Z))\to H^2(X',\Bbb Z/n\Bbb Z)$. Let 
$s'$ be a closed point of $S'$, such that $h(s')=s$. The fibre of the sheaf
$R^1f_{*}(\Bbb Z /n\Bbb Z)$ at $s'$ is isomorphic to $H^1(X_s,\Bbb Z /n\Bbb Z)$. Let $c_s\in
H^1(X_s,\Bbb Z /n\Bbb Z)$ be the class corresponding to the torsor 
$f_s:Y_s\to X_s:=X\times _Sk(s)$. Then $c_s$ can be lifted to a global section $c\in
H^0(S',R^1f'_*(\Bbb Z/n\Bbb Z))$, since $R^1f_{*}(\Bbb Z /n\Bbb Z)/S'$ is constant.
The element $c$ is the image of a class $\tilde c\in H^1(X',\Bbb Z/n\Bbb Z)$, via the above 
sequence, if and only if its image $c'$ in  $H^2(S',F'_*(\Bbb Z/n\Bbb Z))$ vanishes.
The element $c'$ is thus the obstruction to lift the $\mu_n$-torsor 
$f_s:Y_s\to X_s:=X\times _Sk(s)$ to a $\mu_n$-torosr $f':Y'\to X'$. We will show that 
$c'=0$. The image of $c'$ in  $H^2(X',\Bbb Z/n\Bbb Z)$ vanishes, thus it also vanishes
in $H^2(Y',\Bbb Z/n\Bbb Z)$ via the canonical map $H^2(X',\Bbb Z/n\Bbb Z)\to 
H^2(Y',\Bbb Z/n\Bbb Z)$. We have a canonical map 
$H^2(S',\Bbb Z/n\Bbb Z)\to H^2(Y',\Bbb Z/n\Bbb Z)$. We also have a norm 
map $H^2(Y',\Bbb Z/n\Bbb Z)\to H^2(S',\Bbb Z/n\Bbb Z)$, and the composite map 
$H^2(S',\Bbb Z/n\Bbb Z)\to H^2(Y',\Bbb Z/n\Bbb Z)\to H^2(S',\Bbb Z/n\Bbb Z)$                        is multiplication by $d$, Hence, we deduce that the class $c'$ is annihilated by $d$.
Since it is also annihilated by $n$, the group $H^2(S',\Bbb Z/n\Bbb Z)$ being $n$-torsion,
we deduce that $c'=0$.

\pop {.5}
\par
\noindent
{\bf \gr 6.8. Lemma.}\rm\ {\sl Let $k$ be an algebraically 
closed field, of characteristic $p$. Let $S$ be a smooth, complete, and 
irreducible $k$-curve. Let $f:X\to S$ be a smooth family of 
curves, of genus $g\ge 2$. Assume that the fundamental group 
$\pi_1$ is constant on the family $f$. Then for every finite cover 
$S'\to S$, and every finite \'etale cover $Y'\to X':=X\times _SS'$, 
the fundamental group $\pi_1$, is also constant on the family $Y'\to S'$.}

\pop {.5}
\par
\noindent
{\bf \gr Proof.} Standard, using the functorial properties of 
fundamental groups.

\pop {.5}
\par
\noindent
{\bf \gr 6.9. Lemma.}\rm\ {\sl Let $k$ be an algebraically 
closed field, of characteristic $p$. Let $S$ be a smooth, complete, and 
irreducible $k$-curve. Let $f:X\to S$ be a smooth family of 
curves, of genus $g\ge 2$. Let $S'\to S$ be a finite cover, and let 
$Y'\to X':=X\times _SS'$ be an \'etale cover. Assume that the smooth relative
$S'$-curve $Y'\to S'$
 is isotrivial. Then, the smooth relative $S$-curve 
$X\to S$ is also isotrivial}.

\pop {.5}
\par
\noindent
{\bf \gr Proof.}\rm\ Use Lemma 1.32,
 in [Ta-2].

\pop {.5}
\par
Next, we will prove the main Theorem 6.5.

\pop {.5}
\par
\noindent
{\bf \gr Proof of Theorems 6.5 and 6.6.}\rm\ 
Fix a closed point $s$ of $S$, and let 
$X_s:=X\times _Sk(s)$ be the fibre of $X$ above the point $s\in S$. 
By Tamagawa's result, Theorem 3.2,  we can find an 
\'etale cover $Y_s\to X_s$, such that the gonality
$d_{Y_s}$ of $Y_s$ is $\ge 5$. Moreover, 
the cover $Y_s\to X_s$ can be chosen to
be a composition of two cyclic, of order prime to $p$, covers. In particular, 
we can find, by Lemma 6.7, a finite \'etale cover $h:S_1\to S$, and an \'etale cover
$f_1:Y\to X_1:=X\times _S S_1$, such that the pair $(f_1,h)$ 
is a good lifting of the cover $Y_s\to X_s$. Now, consider the 
smooth family of curves $Y\to S_1$. Let $s_1$ be a closed point of $S_1$, 
such that $h(s_1)=s$. Then, by construction, the gonality of the fibre
$Y_{s_1}:=Y\times_S{k(s_1)}$ of $Y$ above $s_1$ is $\ge 5$. In particular, 
by using Theorem 3.3, we can find a prime integer $l>>0$, 
distinct from $p$, and a non trivial $\mu_l$ torsor $Z_{s_1}\to Y_{s_1}$, 
such that the following two conditions hold:
\pop {.5}
\par
i)\ The $\mu_l$-torsor $Z_{s_1}\to Y_{s_1}$ is new ordinary.
\pop {.5}
\par
ii)\ The natural map, $T:M_{Y_{s_1}}\to M_{J_{s_1}^{\new}}$, induced by the 
$\mu_l$-torsor $Z_{s_1}\to Y_{s_1}$, where $J_{s_1}^{\new}$ is the new part
of the jacobian of $Z_{s_1}$, with respect to the above torsor, is 
an immersion (cf. 3).
\pop {.5}
\par
Also, by applying Lemme  6.7, 
we can find a finite \'etale cover $h':S_2\to S_1$, and a $\mu_l$-torsor
$f_2:Z'\to Y':=Y\times _{S_1} S_2$, such that the pair $(f_2,h')$, is 
a good lifting of the $\mu_l$-torsor $Z_{s_1}\to Y_{s_1}$.
Let $J_{Y'}$ (resp. $J_{Z'}$) be the relative jacobian
of $Y'$ over $S_2$ (resp. the relative jacobian of $Z'$ over $S_2$) 
which is an $S_2$-abelian scheme. Let $f_2^*:J_{Y'}\to J_{Z'}$ be the canonical
homomorphism, which is induced by the pull back of invertible sheaves.
Let $J^{\new}:=J_{Z'}/f_2^*(J_{Y'})$, 
be the new part of the jacobian $J_{Z'}$, with respect to the $\mu_l$-torsor 
$Z'\to Y'$.
Let $\eta'$ be the generic point of $S_2$, and let $s_2$ be a point of 
$S_2$ such that $h'(s_2)=s_1$. The fibre $J_{s_2}^{\new}:=
J^{\new}\times _{S_2} k(s_2)$ of $J^{\new}$, 
above the point $s_2$, is by construction an ordinary abelian variety. 
This implies, a fortiori, that the generic fibre $J_{\eta'}^{\new}:=
J^{\new}\times _{S_2} k(\eta')$ of $J^{\new}$, where $\eta'$ is the generic 
point of $S_2$, is also ordinary, since $J_{\eta'}^{\new}$ specializes 
to $J_{s_2}^{\new}$. The following two cases can occur:
\pop {.5}
\par
\ {\bf Case 1}:\ The abelian scheme $J^{\new}\to S_2$ has constant 
$p$-rank, i.e. all fibres of $J^{\new}$ over $S_2$ 
are ordinary abelian varieties. Then,
since $S_2$ is complete, we deduce from Theorem 4.2, that the abelian scheme
$J^{\new}\to S_2$ is isotrivial. Note, that the deformation $J^{\new}$
of $J_{s_2}^{\new}$, induces a first order infinitesimal deformation
of $J_{s_2}^{new}$, which is a trivial deformation since $J^{\new}$ is isotrivial, and 
which by construction is the image of the first order
infinitesimal deformation of $Y_{s_1}$ induced by $Y'$, 
via the above natural map $T:M_{Y_{s_1}}\to M_{J_{s_1}^{\new}}$.
Since, by construction, the map $T$ is an immersion, we conclude that the 
deformation $Y'\to S_2$ is isotrivial, as well. A fortiori, 
the family $X\to S$ is also isotrivial, by Lemma 6.9. But this contradicts
our hypothesis that the family $X\to S$ is not isotrivial. So, 
case 1 can not occur.

\pop {.5}
\par
\ {\bf Case 2}:\ The abelian scheme $J^{\new}\to S_2$ does not have constant 
$p$-rank, i.e. there exists a closed point $\tilde s\in S_2$, such that 
the $p$-rank of the fibre $J_{\tilde s}^{\new}:=J^{\new}\times _{S_2}
\tilde s$ of $J^{\new}$ over the point $\tilde s$, is strictly smaller than
the $p$-rank of the generic fibre $J_{\eta'}^{\new}$ of $J^{\new}$. This, in
particular, implies that the $p$-rank of the fibre $Z'_{\tilde s}:=Z'\times_
{S_2}{\tilde s}$, of $Z'$ above the point $\tilde s$, is strictly smaller than
the $p$-rank of the generic fibre $Z'_{\eta '}:=Z'\times _{S_2}{\eta '}$ 
of $Z'$. This already proves Theorem 6.6. Now, this implies in particular that 
the geometric fundamental group $\pi_1$ is not constant on the family
$Z'\to S_2$. Thus, by Lemma 6.8, we deduce that the geometric fundamental 
group $\pi_1$ is not constant on the family $X\to S$. This finishes 
the proof of Theorem 6.5.

\pop {.5}
\par
Theorem 6.5, can be generalized, to the situation where we consider the full 
tame fundamental group. Note, first, that given a family of curves 
$f:X\to S$ of genus $g\ge 0$, and $n$ sections $\{s_1,s_2,...,s_n\}$, of 
the morphism $f$, we can extend the notion of isotriviality, as in 
Definition 5.2, to the pair $(f,\{s_1,s_2,...,s_n\})$, 
by considering the moduli scheme $\Cal M_{g,n}$, of $n$-pointed smooth 
and projective curves of genus $g$. Now, it is easy to deduce from 6.5 
the following theorem:

\pop {.5}
\par
\noindent
{\bf \gr 6.10. Theorem.}\rm\ {\sl Let $k$ be an algebraically closed field, of 
characteristic $p$. Let $S$ be a smooth, complete, and irreducible $k$-curve,
with generic point $\eta$. Let $f:X\to S$ be a smooth family of curves, 
of genus $g\ge 0$. Let $\{s_1,s_2,...,s_n\}$ be $n$ sections of $f$, 
with disjoint support. Assume that the pair $(f,\{s_1,s_2,...,s_n\})$ is not 
isotrivial. Then, the tame fundamental group is not constant on the pair 
$(f,\{s_1,s_2,...,s_n\})$, i.e. there exists a closed point $s\in S$, such
that the specialization homomorphism $\Sp:\pi_1^t(X_{\eta}-{\{s_1(\eta),
s_2(\eta),...,s_n(\eta)\}})\to \pi_1^t(X_{s}-{\{s_1(s),s_2(s),...,s_n(s)\}})$
between tame fundamental groups, where $X_{\eta}:=X\times _S\eta$ 
(resp.  $X_{s}:=X\times _S s$), is not an isomorphism.}

\pop {1}
\par
\noindent
{\bf \gr References.}

\pop {.5}
\par
\noindent
[Bo-Lu-Ra] S. Bosch W. L\"utkebohmert and M. Raynaud, 
N\'eron Models, Ergebnisse der Mathematik, 3. Folge, Band 21, (1990).

\pop {.5}
\par
\noindent
[Da-Sh] V. I. Danilov and V. V. Shokurov, Algebraic Curves, 
Algebraic Manifolds and Schemes, Springer, (1998).

\pop {.5}
\par
\noindent
[De-Mu] P. Deligne and D. Mumford, {\sl The irreducibility of the space 
of curves with a given genus}, Publ. Math. IHES, 36, 75-110, (1969).

\pop {.5}
\par
\noindent
[Fa-Lo] C. Faber, E. Looijenga, {\sl Remark on moduli of curves}, In Moduli 
of Curves and Abelian Varieties, C. Faber and E. Looijenga (Eds.), 
Aspects of Mathematics, E 33, (1999).

\pop {.5}
\par
\noindent
[Gr] A. Grothendieck, Brief an G. Faltings (letter to G. Faltings), Geometric 
Galois Action, 1, 49-58, London Math. Soc. Lecture Notes ser., 
Cambridge Univ. Press, Cambridge (1997).

\pop {.5}
\par
\noindent
[Mi] J.S. Milne, {\sl Jacobian Varieties}, In Arithmetic Geometry, G. Cornell
and J. H. Silverman (Eds), Springer-Verlag, (1985).

\pop {.5}
\par
\noindent
[Mo] L. Moret-Bailly, Pinceaux de vari\'et\'es ab\'eliennes, Ast\'erisque, 129,
(1985).

\pop {.5}
\par
\noindent
[Mu] D. Mumford, Curves and their Jacobians, Ann Arbor, The University 
Of Michigan Press, (1976)

\pop {.5}
\par
\noindent
[No-Oo] P. Norman and F. Oort, {\sl Moduli of abelian varieties}, 
Annals of Math. 112, 413-439, (1980)

\pop {.5}
\par
\noindent
[Oo] F. Oort, {\sl A Stratification of a moduli space of polarized abelian 
varieties in positive characteristic}, In Moduli of Curves and 
Abelian Varieties, C. Faber and E. Looijenga (Eds.), 
Aspects of Mathematics E 33, (1999).

\pop {.5}
\par
\noindent
[Oo-1] F. Oort, {\sl Subvarieties of moduli spaces}, Inventiones math., 24,
95-119 (1974).

\pop {.5}
\par
\noindent
[Po-Sa] F. Pop and M. Sa\"\i di, {\sl On the specialization homomorphism
of fundamental groups of curves in positive characteristic}, To appear
in the Proceedings of the MSRI conference, on Arithmetic fundamental groups, 
(1999)

\pop {.5}
\par
\noindent
[Pr] R. Pries, {\sl Families of wildly ramified covers of curves}, Amer. 
J. Math, Vol. 124, 4, 737-768, (2002).

\pop {.5}
\par
\noindent
[Ra] M. Raynaud, {\sl Section des fibr\'es vectoriel sur une courbe}, 
Bull. Soc. Math. France, t. 110, 103-125 (1982).

\pop {.5}
\par
\noindent
[Ra-1] M. Raynaud, {\sl Sur le groups fondamental d'une courbe compl\`ete
en caract\'eristique $p>0$}, To appear in the proceedings, 
of the MSRI conference, on arithmetic fundamental groups (1999).

\pop {.5}
\par
\noindent
[SGA-1] A. Grothendieck, Rev\^etements \'etales et groupe fondamental, 
L.N.M. Vol. 224, Springer Verlag, (1971).

\pop {.5}
\par
\noindent
[Se] J. P. Serre, {\sl Sur la topologie des vari\'et\'es alg\'ebriques en 
caract\'eriqstique $p>0$}, Oeuvres Compl\'etes Vol. 1. 

\pop {.5}
\par
\noindent
[Se-1] J. P. Serre, Algebraic groups and class fields, Graduate Texts in 
Mathematics, Springer Verlag, (1975).

\pop {.5}
\par
\noindent
[Sh] I. Shafarevich, {\sl On $p$-extensions}, A.M.S. translation, 
Serie 2, Vol. 4 (1965).

\pop {.5}
\par
\noindent
[Sz] L. Szpiro, {\sl Propri\'et\'es num\'eriques du faisceau dualisant 
relatif}, In S\'eminaire sur les pinceaux de courbes de genre au moins deux,
Ast\'erisque 86, (1981).

\pop {.5}
\par
\noindent
[Ta] A. Tamagawa, {\sl On the tame fundamental groups of curves over 
algebraically closed fields of characteristic $>0$}, preprint.

\pop {.5}
\par
\noindent
[Ta-1] A. Tamagawa, {\sl Finiteness of isomorphism classes of curves in 
positive characteristic with prescribed fundamental group}, preprint.

\pop {.5}
\par
\noindent
[Ta-2] A. Tamagawa, {\sl Fundamental groups and geometry of curves 
in positive characteristic}, To appear in Proceedings of Symposia in Pure
Mathematics.

\pop {.5}
\par
\noindent
[Ta-3] A. Tamagawa, {\sl The Grothendieck conjecture for affine curves},
Compositio Math., 109, 135-194, (1997)

\pop {2}
Mohamed Sa\"\i di

\pop {1}
\par
Departement of Mathematics
\par
University of Durham
\par
Science Laboratories
\par
South Road
\par
Durham DH1 3LE United Kingdom
\par
Mohamed.Saidi\@durham.ac.uk

\end
\enddocument